%
%
%
%
\magnification=\magstep1
\input amstex
\documentstyle{amsppt}
\voffset-3pc

\define\C{{\Bbb C}}

\define\dee{\partial}
\redefine\O{\Omega}
\redefine\phi{\varphi}
\define\Obar{\overline{\Omega}}
\define\Ot{\widetilde\Omega}
\define\Oh{\widehat\Omega}

\NoRunningHeads
\topmatter
\title
The Bergman kernel and quadrature domains in the plane
\endtitle
\author Steven R. Bell${}^*$ \endauthor
\thanks ${}^*$Research supported by NSF grant DMS-0305958 \endthanks
\keywords Bergman kernel, Szeg\H o kernel
\endkeywords
\subjclass 30C40 \endsubjclass
\address
Mathematics Department, Purdue University, West Lafayette, IN  47907 USA
\endaddress
\email bell\@math.purdue.edu \endemail
\abstract
A streamlined proof that the Bergman kernel associated to a
quadrature domain in the plane must be algebraic will be given.
A byproduct of the proof will be that the Bergman kernel is a
rational function of $z$ and one other explicit function known
as the Schwarz function.  Simplified proofs of several other
well known facts about quadrature domains will fall out along
the way.  Finally, Bergman representative coordinates will be
defined that make subtle alterations to a domain to convert it
to a quadrature domain.  In such coordinates, biholomorphic
mappings become algebraic.
\endabstract
\endtopmatter
\document

\hyphenation{bi-hol-o-mor-phic}
\hyphenation{hol-o-mor-phic}

\subhead 1. Introduction \endsubhead
In this paper, we will recombine a string of results by Aharonov and
Shapiro \cite{1}, Gustafsson \cite{14}, Davis \cite{12}, Shapiro
\cite{17}, and Avci \cite{2} in light of recent results by the author
in \cite{7} and \cite{9} to obtain elementary proofs of a number of
results about quadrature domains and the classical functions associated
to them.  In particular, we present an efficient proof of the fact proved
in \cite{9} that the Bergman kernel associated to a quadrature
domain in the plane is an algebraic function.  In fact, we shall
show that it is a rational function of $z$ and the Schwarz function,
and consequently it is also a rational combination of $z$ and $Q(z)$
where $Q(z)$ is an explicit algebraic function given by
$$Q(z)=\int_{b\O}\frac{\bar w}{w-z}\ dw.$$
These results are all a natural outgrowth of the work of Aharonov
and Shapiro \cite{1} and Gustafsson \cite{14} and
many of the results obtained in those works will come as corollaries in
the approach we take here.  For example, Aharonov and Shapiro proved
that Ahlfors maps associated to quadrature domains are algebraic,
and we shall deduce this via the connection between Ahlfors maps,
the double, and the Bergman kernel.  The new approach we use shall also
allow us to view Bergman representative coordinates in a new and
very interesting light.

We shall call an $n$-connected domain
$\O$ in the plane such that no boundary component is a point a {\it
quadrature domain\/} if there exist finitely many points
$\{w_j\}_{j=1}^N$ in the domain and non-negative integers $n_j$ such
that complex numbers $c_{jk}$ exist satisfying
$$\int_\O f\ dA = \sum_{j=1}^N\sum_{k=0}^{n_j} c_{jk} f^{(k)}(w_j)\tag1.1$$
for every function $f$ in the Bergman space of square integrable
holomorphic functions on $\O$.  Here, $dA$ denotes Lebesgue area
measure.  Our results require the function $h(z)\equiv1$ to be
in the Bergman space, and so we shall also assume that the domain
under study has {\it finite area}.

If $\O$ is an $n$-connected quadrature domain of finite area in the
plane such that no boundary component is a point, then it is well
known that the Bergman kernel function associated to $\O$ satisfies
an identity of the form
$$1\equiv\sum_{j=1}^N\sum_{m=0}^{n_j} c_{jm} K^{(m)}(z,w_j)
\tag1.2$$
where $K^{(m)}(z,w)$ denotes $(\dee^m/\dee\bar w^m)K(z,w)$ (and of
course $K^{(0)}(z,w)=K(z,w)$) and where the points $w_j$ are the points
that appear in the characterizing formula (1.1) of quadrature domains.
It can be seen by noting that the inner product of an analytic
function against the function $h(z)\equiv1$ and against the sum on
the right hand side of (1.2) agree for all functions in the Bergman
space.  Hence the two functions must be equal.  Note that we must
assume that $\O$ has finite area here just so that $h(z)\equiv1$ is
in the Bergman space.

We first state a theorem about the Bergman kernel of a multiply
connected domain in the plane with smooth boundary.  Note that,
although formula (1.2) is clearly in the background, we do not
assume that the domain is a quadrature domain.

\proclaim{Theorem 1.1}
Suppose that $\O$ is an $n$-connected bounded domain in the plane
whose boundary is given by $n$ non-intersecting simple closed
$C^\infty$ smooth real analytic curves.  Let $A(z)$ be a function
of the form
$$\sum_{j=1}^N\sum_{m=0}^{n_j} c_{jm} K^{(m)}(z,w_j)$$
where the $w_j$ are points in $\O$.  Let $G_1$ and $G_2$ be any two
meromorphic functions on $\O$ that extend meromorphically to the
double of $\O$ and form a primitive pair for the field of
meromorphic functions on the double.  There is a rational
function $R(z_1,z_2,w_1,w_2)$ of four complex variables
such that $K(z,w)$ is given by
$$K(z,w)=A(z)\overline{A(w)}R\left(G_1(z),G_2(z),
\overline{G_1(w)},\overline{G_2(w)}\right).$$
\endproclaim

Theorem~1.1 is a corollary of Theorem~2.3 of \cite{7}.  We
shall give a straightforward and direct proof of Theorem~1.1
in \S3.

We remark here that two meromorphic functions on a compact
Riemann surface are said to form a primitive pair if they
generate the field of meromorphic functions, i.e., if every
meromorphic function on the Riemann surface is a rational
combination of the two.  For basic facts about primitive
pairs, see Farkas and Kra \cite{13}.

If $\O$ is an $n$-connected domain in the plane such that no
boundary component is a point, then it is a standard construction
in the subject to produce a biholomorphic mapping $\phi$ which
maps $\O$ one-to-one onto a bounded domain $\O_s$ bounded by $n$
smooth real analytic curves.  The subscript $s$ stands for ``smooth''
and we shall write $K(z,w)$ for the Bergman kernel of $\O$ and
$K_s(z,w)$ for the Bergman kernel associated to $\O_s$.  We shall
say that a meromorphic function $h$ on $\O$ {\it extends meromorphically
to the double} of $\O$ if $h\circ\phi^{-1}$ is a meromorphic function
on $\O_s$ which extends meromorphically to the double of $\O_s$.
(This terminology might be considered to be rather non-standard,
but it greatly simplifies the statements of many of our results below.)
It is easy to verify that this definition does not depend on the
choice of $\phi$ and $\O_s$.  We shall also say that $G_1$ and
$G_2$ form a primitive pair for $\O$ if $G_1\circ\phi^{-1}$ and
$G_2\circ\phi^{-1}$ extend meromorphically to the double of $\O_s$
and form a primitive pair for the double of $\O_s$.  Using this
terminology, we may apply Theorem~1.1 to obtain the following
result.

\proclaim{Theorem 1.2}
Suppose that $\O$ is an $n$-connected quadrature domain of finite area
in the plane such that no boundary component is a point.  Then the
Bergman kernel function $K(z,w)$ associated to $\O$ is a rational
combination of any two functions $G_1$ and $G_2$ that form a primitive
pair for $\O$ in the sense that $K(z,w)$ is a rational combination of
$G_1(z)$, $G_2(z)$, $\overline{G_1(w)}$, and $\overline{G_2(w)}$.
\endproclaim

We remark that the same conclusion in Theorem~1.2 can be made about
the square $S(z,w)^2$ of the Szeg\H o kernel.  Furthermore, the
classical functions $F_j'$ (see \S2 for definitions) are rational functions
of $G_1$ and $G_2$ and so is any proper holomorphic mapping of $\O$
onto the unit disc.  The reader may see \cite{7} for the details of
the more general statement.

The next theorem of Gustafsson \cite{14} reveals that quadrature
domains have nice boundaries.  It also allows us to smoothly connect
a quadrature domain to another domain with a double in the classical
sense.

\proclaim{Theorem 1.3}
Suppose that $\O$ is an $n$-connected quadrature domain of finite area
in the plane such that no boundary component is a point.  Suppose
that $\phi$ is a holomorphic mapping which maps $\O$ one-to-one onto
a bounded domain $\O_s$ bounded by $n$ smooth real analytic curves.
Then $\phi^{-1}$ extends holomorphically past the boundary of $\O_s$.
Furthermore, $\phi^{-1}$ extends meromorphically to the double of
$\O_s$.  It follows that the boundary of $\O$ is piecewise real analytic
and the possibly finitely many non-smooth boundary points are easily
described as cusps which point toward the inside of $\O$.  Furthermore,
it follows that $\phi$ extends continuously to the boundary of $\O$.
\endproclaim

Avci did not state Theorem~1.3 in \cite{2}, however, all the
elements of a proof are there.  If Avci had done things in a different
order, he might have needed this theorem and he might very well have
spelled out the proof as well.  In fact, much of Avci's work in
\cite{2} is headed in a direction that could have easily have led
to many of the results of this paper.

We shall give a short Bergman kernel proof of Theorem~1.3 in the \S3.

We remark that at a cusp boundary point $b$ of $\O$ mentioned in
Theorem~1.3, the map $\phi$ behaves like a principal branch of the
square root of $z-b$ mapping the plane minus a horizontal slit from
$b$ to the left to the right half plane.  This map sends $b$ to zero,
and the inverse map $\phi^{-1}$ behaves like $b+z^2$ on the right
half plane mapping zero to $b$.  The power two in $b+z^2$ is the
only power larger than one that makes such a map one-to-one on
the right half plane near zero, and this kind of reasoning can be
used to show that the derivative of $\phi^{-1}$ can have at most a
simple zero at a boundary point of $\O_s$.

Aharonov and Shapiro \cite{1} showed that the Schwarz function
associated to a quadrature domain extends meromorphically to the
domain, i.e., that the function $\bar z$ agrees on the boundary with
a function $S(z)$, known as the Schwarz function, which is meromorphic
on the domain and which extends continuously up to the boundary.  We
will follow Gustafsson \cite{14} and modify somewhat his observation
that $z$ and $S(z)$ form a primitive pair (in the special sense we use
here) to obtain a quick proof of the following result.

\proclaim{Theorem 1.4}
Suppose that $\O$ is an $n$-connected quadrature domain of finite area
in the plane such that no boundary component is a point.  The
function $z$ extends to the double as a meromorphic function.
Consequently, there is a meromorphic function $S(z)$ on $\O$ (known as
the Schwarz function) which extends continuously up to the boundary
such that $S(z)=\bar z$ on $b\O$.  The functions $z$ and $S(z)$ form a
primitive pair for $\O$.    It also follows that $z$ and
$$Q(z)=\int_{b\O}\frac{\bar w}{w-z}\ dw$$
form a primitive pair for $\O$.  Consequently, $S(z)$ and $Q(z)$ are
algebraic functions.  It also follows that the boundary of $\O$ is a
real algebraic curve.
\endproclaim

Aharonov and Shapiro \cite{1} first showed that the boundary of a
quadrature domain is an algebraic curve and Gustafsson \cite{14}
later gave a precise description of what these curves must be.
Theorems~1.4 and~1.2 can now be combined to yield the following
result.

\proclaim{Theorem 1.5}
Suppose that $\O$ is an $n$-connected quadrature domain of finite area
in the plane such that no boundary component is a point.  Then the
Bergman kernel function $K(z,w)$ associated to $\O$ is a rational
combination of the two functions $z$ and the Schwarz function.  It is
also a rational function of $z$ and $Q(z)$.
Consequently, $K(z,w)$ is algebraic.  Furthermore, the Szeg\H o kernel
is algebraic, the classical functions $F_j'$ are algebraic, and every
proper holomorphic mapping from $\O$ onto the unit disc is algebraic.
\endproclaim

Similar statements to the theorems above can be made for the Poisson
kernel and first derivative of the Green's function.  These results
follow from formulas appearing in \cite{6} and we do not spell them
out here.

It is interesting to note here that, not only is the Bergman kernel
of a quadrature domain $\O$ a rational combination of $z$ and the
Schwarz function, but the Schwarz function associated to $\O$ is a
rational combination of $K(z,a)$ and $K(z,b)$ for two fixed points
$a$ and $b$ in the domain by virtue of the fact that $S(z)$ extends
to the double and because it is possible to find two such functions
$K(z,a)$ and $K(z,b)$ of $z$ that form a primitive pair for $\O$.

Since $S(z)=\bar z$ on the boundary, the next theorem is an easy
consequence of Theorem~1.5.

\proclaim{Theorem 1.6}
Suppose that $\O$ is an $n$-connected quadrature domain in the plane
of finite area such that no boundary component is a point. The Bergman
kernel $K(z,w)$ and the square $S(z,w)^2$ of the Szeg\H o kernel are
rational functions of $z$, $\bar z$, $w$, and $\bar w$ on
$b\O\times b\O$ minus the boundary diagonal.  The functions $F_j(z)$
are rational functions of $z$ and $\bar z$ when restricted to the boundary.
Furthermore, the unit tangent vector function $T(z)$ is such that
$T(z)^2$ is a rational function of $z$ and $\bar z$ for $z\in b\O$.
\endproclaim

The Riemann mapping theorem can be viewed as saying that any
simply connected domain in the plane that is not the whole plane
is biholomorphic to the grandaddy of all quadrature domains, the
unit disc.  Gustafsson generalized this theorem to multiply connected
domains.  He proved that any finitely connected domain in the plane
such that no boundary component is a point is biholomorphic to a
quadrature domain.  The circle of ideas we develop in this paper
can be used to reformulate these theorems via a tool I would venture
to call ``Bergman representative coordinates.''  We shall show that
any domain with $C^\infty$ smooth boundary can be mapped by a
biholomorphic mapping which is as $C^\infty$ close to the identity
map as desired to a quadrature domain.  The biholomorphic map
will be given in the form of a quotient of linear combinations
of the Bergman kernel.  This process begins with the following
lemma, which was proved in \cite{3}.  Let $A^\infty(\O)$ denote the
subspace of $C^\infty(\Obar)$ consisting of functions that are
holomorphic on $\O$.

\proclaim{Lemma 1.7}
Suppose that $\O$ is a bounded finitely connected
domain bounded by simple closed $C^\infty$ smooth curves.  The
complex linear span of the set of functions of $z$ of the form
$K(z,b)$ where $b$ are points in $\O$ is dense in $A^\infty(\O)$.
\endproclaim

If $\O$ is a domain as in Lemma~1.7, let $\Cal K_2(z)$ denote a finite
linear combination $\sum c_jK(z,b_j)$ which is $C^\infty$ close to the
function $h(z)\equiv 1$ and let
$\Cal K_1(z)$ denote a finite linear combination $\sum a_jK(z,b_j)$
which is $C^\infty$ close to the
function $h(z)\equiv z$.  I call the quotient
$\Cal K_1(z)/\Cal K_2(z)$ a Bergman representative mapping function if
it is one-to-one.  By taking the linear combinations to be $C^\infty$
close enough to their target functions, any such mapping can be made
one-to-one and as $C^\infty$ close to the indentity as desired.
The mappings extend meromorphically to the double of $\O$ because
$T(z)K(z,b)=-\overline{T(z)\Lambda(z,b)}$ for $z$ in $b\O$.
Hence, according to Gustafsson \cite{14}, the mapping sends $\O$
to a quadrature domain.  Hence, the following theorem holds

\proclaim{Theorem 1.8}
Suppose that $\O$ is a bounded finitely connected
domain bounded by simple closed $C^\infty$ smooth curves.  There
is a Bergman representative mapping function which is as $C^\infty$
close to the identity map as desired which maps $\O$ to a quadrature
domain.
\endproclaim

Bergman representative mappings as defined here yield a rather
fascinating change of coordinates.  Indeed, they were used in
several complex variables in \cite{10} to locally linearize
biholomorphic mappings.  In one variable, if $\Phi:\O_1\to\O_2$ is
a biholomorphic (or merely proper holomorphic mapping) between
$C^\infty$-smooth finitely connected domains in the plane, then
Theorem~1.8 allows us to make changes of coordinates that are
$C^\infty$ close to the identity on each domain in such a way
that the mapping $\Phi$ in the new coordinates is an {\it algebraic
function}.  Indeed, if $\Phi:\O_1\to\O_2$ is a proper holomorphic
mapping between quadrature domains in the new coordinates, let
$f_a$ denote an Ahlfors mapping of $\O_2$ onto the unit disc.
Then $f_a$ is a proper holomorphic mapping of the smooth quadrature
domain $\O_2$ onto the unit disc, and is therefore algebraic by
Theorem~1.5.  Now $f_a\circ\Phi$ is a proper holomorphic mapping
of the smooth quadrature domain $\O_1$ onto the unit disc, and is
therefore algebraic.  It follows that $\Phi$ itself must be
algebraic.  It is rather striking that subtle changes in the boundary
can make conformal mappings become defined on the whole complex plane.

Since Bergman representative mappings extend to the double and
are one-to-one, they are Gustafsson mappings (as defined in
\cite{9}), and they can be used to compress the classical
kernel functions into a small data set as in \cite{9}.
Bj\"orn Gustafsson read a preliminary version of this
paper and realized that every Gustafsson map can be expressed
as a Bergman representative map.  He has granted me permission to
include his argument in \S4 of the present paper.

The points in the quadrature identity associated to the Bergman
representative domain in Theorem~1.8 can be arranged to fall in
any small disc in the domain.  This can be done using similar
constructions to those used in \cite{9} and by Gustafsson in
\cite{14}.  We do not treat this problem here.

Another interesting way to view Theorem~1.8 is as follows.  Shrink
a smooth domain $\O$ by moving in along an inward pointing unit
vector a fixed short distance.  Now use a Bergman
representative mapping which is sufficiently $C^\infty$ close to the
identity map so that the shrunken domain gets mapped to a domain
inside of $\O$.  This shows that we may lightly ``sand'' the edges
of our original domain to turn it into a quadrature domain.
Similarly, by expanding the domain by first moving along an outward
pointing normal and repeating this process, we can see that we
can ``paint'' the edges of our domain with an arbitrarily thin
coat of paint of variable thickness to turn it into a quadrature
domain.

It is reasonable to allow more general functions of the form
$$\sum_{j=1}^N\sum_{m=0}^{n_j} c_{jm} K^{(m)}(z,w_j)$$
to appear as the numerator and denominator of a Bergman representative
mapping because the individual elements of these functions
satisfy the same kind of basic identity as the Bergman kernel
itself (see formula (3.1)) and, consequently, such
quotients extend to the double and have the same mapping properties
as the Bergman representative mappings as we constructed them above.
Under this stipulation, we can also prove a converse to Theorem~1.8.

\proclaim{Theorem 1.9}
If $\O$ is a finitely connected quadrature domain of finite
area, then there is a Bergman representative mapping which is
equal to the identity map.
\endproclaim

Another approach to representative coordinates is given by Jeong
and Taniguchi in \cite{15}.  It might be interesting to see what
comes of combining the two approaches.

Finally, we remark that, when the results of this paper are combined
with the results in \cite{8}, it can be seen that the infinitesimal
Carath\'eodory metric associated to a finitely connected quadrature
domain of finite area such that no boundary component is a point is
given by $\rho(z)|dz|$ where $\rho$ is a rational combination of $z$,
the Schwarz function, and the complex conjugates of these two functions.

Complete proofs of the theorems will be given in \S3.

\subhead 2. Preliminaries\endsubhead
It is a standard construction in the theory of conformal mapping to
show that an $n$-connected domain $\O$ in the plane such that no boundary
component is a point is conformally equivalent via a map $\phi$ to a
bounded domain $\Ot$ whose boundary consists of $n$ simple closed
$C^\infty$ smooth real analytic curves.  Since such
a domain $\Ot$ is a bordered Riemann surface, the double of $\Ot$
is an easily realized compact Riemann surface.  We shall say that
an analytic or meromorphic function $h$ on $\O$ {\it extends meromorphically
to the double of\/} $\O$ if $h\circ\phi^{-1}$ extends meromorphically
to the double of $\Ot$.  Notice that whenever $\O$ is itself a bordered
Riemann surface, this notion is the same as the notion that $h$
extends meromorphically to the double of $\O$.  We shall say that
two functions $G_1$ and $G_2$ extend to the double and generate the
meromorphic functions on the double of $\O$, and that they therefore
form a primitive pair for the double of $\O$, if $G_1\circ\phi^{-1}$ and
$G_2\circ\phi^{-1}$ extend to the double of $\Ot$ and form a primitive
pair for the double of $\Ot$ (see Farkas and Kra \cite{13} for the
definition and basic facts about primitive pairs).

It is proved in \cite{6} that if $\O$ is an $n$-connected domain
in the plane such that no boundary component is a point, then almost
any two distinct Ahlfors maps $f_a$ and $f_b$ generate the meromorphic
functions on the double of $\O$.  It is also proved that any proper
holomorphic mapping from $\O$ to the unit disc extends to the double
of $\O$.

Suppose that $\O$ is a bounded $n$-connected domain whose
boundary consists of $n$ non-intersecting $C^\infty$ smooth simple
closed curves.  The Bergman kernel $K(z,w)$ associated to $\O$ is
related to the Szeg\H o kernel via the identity
$$K(z,w)=4\pi S(z,w)^2+\sum_{i,j=1}^{n-1}
A_{ij}F_i'(z)\overline{F_j'(w)},\tag 2.1$$
where the functions $F_i'(z)$ are classical functions of potential theory
described as follows.  The harmonic function $\omega_j$ which solves the
Dirichlet problem on $\O$ with boundary data equal to one on the boundary
curve $\gamma_j$ and zero on $\gamma_k$ if $k\ne j$ has a multivalued
harmonic conjugate.  Let $\gamma_n$ denote the outer boundary curve.
The function $F_j'(z)$ is a single valued holomorphic function on $\O$
which is locally defined as the derivative of $\omega_j+iv$ where $v$ is
a local harmonic conjugate for $\omega_j$.  The Cauchy-Riemann equations
reveal that $F_j'(z)=2(\dee\omega_j/\dee z)$.

The Bergman and Szeg\H o kernels are holomorphic in the first variable
and antiholomorphic in the second on $\O\times\O$ and they are hermitian,
i.e.,  $K(w,z)=\overline{K(z,w)}$.  Furthermore, the Bergman and
Szeg\H o kernels are in
$C^\infty((\Obar\times\Obar)-\{(z,z):z\in b\O\})$ as functions of $(z,w)$
(see \cite{4, p.~100}).

We shall also need to use the
Garabedian kernel $L(z,w)$, which is related to the Szeg\H o
kernel via the identity
$$\frac{1}{i} L(z,a)T(z)=S(a,z)\qquad\text{for $z\in b\O$ and $a\in\O$}
\tag2.2$$
where $T(z)$ represents the complex unit tangent vector at $z$ pointing
in the direction of the standard orientation of $b\O$.
For fixed $a\in\O$, the kernel $L(z,a)$ is a holomorphic function of $z$
on $\O-\{a\}$ with a simple pole at $a$ with residue $1/(2\pi)$.
Furthermore, as a function of $z$, $L(z,a)$ extends to the boundary
and is in the space $C^\infty(\Obar-\{a\})$.  In fact, $L(z,w)$
is in $C^\infty((\Obar\times\Obar)-\{(z,z):z\in\Obar\})$ as a function
of $(z,w)$ (see \cite{4, p.~102}).  Also, $L(z,a)$ is non-zero for all
$(z,a)$ in $\Obar\times\O$ with $z\ne a$ and $L(a,z)=-L(z,a)$ (see
\cite{4, p.~49}).

For each point $a\in\O$, the function of $z$ given by
$S(z,a)$ has exactly $(n-1)$ zeroes in $\O$ (counting multiplicities) and
does not vanish at any points $z$ in the boundary of $\O$ (see
\cite{4, p.~49}).

Given a point $a\in\O$, the Ahlfors map $f_a$ associated to the pair $(\O,a)$
is a proper holomorphic mapping of $\O$ onto the unit disc.  It is an
$n$-to-one mapping (counting multiplicities), it extends to be in
$C^\infty(\Obar)$, and it maps each boundary curve $\gamma_j$ one-to-one
onto the unit circle.  Furthermore, $f_a(a)=0$, and $f_a$ is the unique
function mapping $\O$ into the unit disc maximizing the quantity $|f_a'(a)|$
with $f_a'(a)>0$.  The Ahlfors map is related to the Szeg\H o kernel
and Garabedian kernel via (see \cite{4, p.~49})
$$f_a(z)=\frac{S(z,a)}{L(z,a)}.\tag2.3$$
Note that $f_a'(a)=2\pi S(a,a)\ne 0$.  Because $f_a$ is $n$-to-one, $f_a$
has $n$ zeroes.  The simple pole of $L(z,a)$ at $a$ accounts for the simple
zero of $f_a$ at $a$.   The other $n-1$ zeroes of $f_a$ are given by the
$(n-1)$ zeroes of $S(z,a)$ in $\O-\{a\}$.

When $\O$ does not have smooth boundary, we define the kernels and domain
functions above as in \cite{5} via a conformal mapping to a domain with
real analytic boundary curves.

\subhead 3. Proofs of the theorems\endsubhead
In this section, we give complete proofs of Theorems~1.1--1.9.

\demo{Proof of Theorem~1.1}
The Bergman kernel is related to the classical Green's function via
(\cite{11, p.~62}, see also \cite{4, p.~131})
$$K(z,w)=-\frac{2}{\pi}\frac{\dee^2 G(z,w)}{\dee z\dee\bar w}.$$
Another kernel function on $\O\times\O$ that we shall need is given by
$$\Lambda(z,w)=-\frac{2}{\pi}\frac{\dee^2 G(z,w)}{\dee z\dee w}.$$
(In the literature, this function is sometimes written as $L(z,w)$
with anywhere between zero and three tildes and/or hats over the top.
We have chosen the symbol $\Lambda$ here to avoid
confusion with our notation for the Garabedian kernel above.)

The Bergman kernel and the kernel $\Lambda(z,w)$
satisfy an identity analogous to (2.2):
$$\Lambda(w,z)T(z)=-K(w,z)\overline{T(z)}\qquad\text{for $w\in\O$ and
$z\in b\O$}\tag3.1$$
(see \cite{4, p.~135}).  We remark that it follows from well known
properties of the Green's function that
$\Lambda(z,w)$ is holomorphic in $z$ and $w$  and is in
$C^\infty(\Obar\times\Obar-\{(z,z):z\in\Obar\})$.
If $a\in\O$, then $\Lambda(z,a)$ has a double pole at
$z=a$ as a function of $z$ and $\Lambda(z,a)=\Lambda(a,z)$ (see
\cite{4, p.~134}).  Since $\O$ has real analytic boundary,
the kernels $K(z,w)$, $\Lambda(z,w)$, $S(z,w)$, and $L(z,w)$,
extend meromorphically to $\Obar\times\Obar$ (see
\cite{4, p.~103, 132--136}).
Let $A(z)=\sum_{j=1}^N\sum_{m=0}^{n_j} c_{jm} K^{(m)}(z,w_j)$
where the $w_j$ are points in $\O$.
Notice that $A$ cannot be the zero function because, if it were, it
would be orthogonal to all functions in the Bergman space, and
consequently every function in the Bergman space would have to
satisfy the identity
$0\equiv\sum_{j=1}^N\sum_{m=0}^{n_j} \bar c_{jm} g^{(m)}(w_j)$, which
is absurd.  Notice that (3.1) shows that there is a meromorphic
function $M(z)$ on $\O$ which extends meromorphically to a
neighborhood of $\Obar$ and which has no poles on $b\O$ such that
$$A(z)T(z)=\overline{M(z)T(z)}\qquad\text{for }z\in b\O.\tag 3.2$$
Let $\Cal B$ denote the class of holomorphic functions $B(z)$ on $\O$
that have the property that they extend holomorphically past the
boundary and such that there exists a meromorphic
function $m(z)$ on $\O$ which extends meromorphically to a
neighborhood of $\Obar$ and which has no poles on $b\O$ such that
$$B(z)T(z)=\overline{m(z)T(z)}\qquad\text{for }z\in b\O.$$
We have shown that $A(z)$ belongs to $\Cal B$.  Notice that, if $B(z)$
is in $\Cal B$, then
$B(z)/A(z)$ is equal to the complex conjugate of a meromorphic
function $m(z)/M(z)$ for $z\in b\O$.  This shows that $B(z)/A(z)$
extends meromorphically to the double of $\O$.
Now the Bergman kernel is given by
$$K(z,w)=4\pi S(z,w)^2+\sum_{i,j=1}^{n-1}
\lambda_{ij}F_i'(z)\overline{F_j'(w)}$$
where
$$S(z,w)=\frac{1}{1-f_a(z)\overline{f_a(w)}}\left(c_0
S(z,a)\overline{S(w,a)}+
\sum_{i,j=1}^{n-1} c_{ij}S(z,a_i)\,\overline{S(w,a_j)}\right)\tag3.3$$
and where $f_a(z)$ denotes the Ahlfors map associated to $a$ (see
\cite{5}).  The functions $F_j'$ belong to the class $\Cal B$.
Indeed
$$F_j'(z)T(z)=-\overline{F_j'(z)T(z)}\qquad\text{for }z\in b\O\tag 3.4$$
(see \cite{4, p.~80}).
Furthermore, functions of $z$ of the form $S(z,a_i)S(z,a_j)$ also
belong to the class $\Cal B$ by virtue of identity (2.2).  Hence,
$K(z,w)$ is given by a sum of terms of the form
$$\frac{B_1(z)\overline{B_2(w)}}{(1-f_a(z)\overline{f_a(w)})^2}$$
plus a sum of functions of the form $B_1(z)\overline{B_2(w)}$
where $B_1$ and $B_2$ belong to $\Cal B$.  Therefore, if we divide
$K(z,w)$  by $A(z)\overline{A(w)}$, we obtain a function which is
a sum of terms of the form
$$\frac{g_1(z)\overline{g_2(w)}}{1-f_a(z)\overline{f_a(w)})^2}$$
plus a sum of functions of the form $g_1(z)\overline{g_2(w)}$
where $g_1$ and $g_2$ extend meromorphically to the double.
But $f_a$ also extends meromorphically to the double.  Hence,
$K(z,w)$ is equal to $A(z)\overline{A(w)}$ times a rational
function of $G_1(z)$, $G_2(z)$, $\overline{G_1(w)}$, and
$\overline{G_2(w)}$ where $G_1$ and $G_2$ are any two functions
that form a primitive pair for the double of $\O$.  This
completes the proof of Theorem~1.1.
\enddemo

\demo{Proof of Theorems~1.2 and~1.3}
Suppose that $\O$ is an $n$-connected quadrature domain of finite area
in the plane such that no boundary component is a point.  Suppose
that $\phi$ is a holomorphic mapping which maps $\O$ one-to-one onto
a bounded domain $\O_s$ bounded by $n$ smooth real analytic curves.
Let $\Phi$ denote $\phi^{-1}$.
The transformation formula for the Bergman kernels under
$\phi$ can be written in the form
$$\Phi'(z)K(\Phi(z),w)=K_s(z,\phi(w))\overline{\phi'(w)}.$$
As mentioned in \S1, since $\O$ is a quadrature domain of finite
area, an identity of the form
$$1\equiv\sum_{j=1}^N\sum_{m=0}^{n_j} c_{jm} K^{(m)}(z,w_j)$$
holds.  Let $A(z)$ denote the linear combination on the right hand side
of this equation.  It now follows that
$$\gather
\Phi'(z)=\Phi'(z)\cdot(A\circ\Phi)(z)=
\Phi'(z)\sum_{j=1}^N\sum_{m=0}^{n_j} c_{jm} K^{(m)}(\Phi(z),w_j) \\
=\sum_{j=1}^N\sum_{m=0}^{n_j} c_{jm} \frac{\dee^m}{\dee\bar w^m}\left(
K_s(z,\phi(w))\overline{\phi'(w)}\right)\left.\vphantom{\int}\right|_{w=w_j}.
\endgather$$
Notice that the function on the last line of this string of equations
is a finite linear combination of functions of the form $K_s^{(m)}(z,W)$
where the $W$ are points in $\O_s$.  Let $A_s(z)$ denote this linear
combination.  We have shown two things.  We have shown that
$$\Phi'=A_s$$
where $A_s$ is a linear combination of functions of $z$ of the form
$K_s^{(m)}(z,W)$.  We shall use this first fact momentarily to prove
Theorem~1.3.  Secondly, we have shown that
$$\Phi'(z)\cdot(A\circ\Phi)(z)=A_s(z),$$
and consequently, that
$$A(z)=\phi'(z)\cdot(A_s\circ\phi)(z).$$
This second fact will yield Theorem~1.2.  Indeed, we may now combine
this fact with Theorem~1.1 and the transformation formula for
the Bergman kernels to obtain that
$$\gather
K(z,w)= \frac{K(z,w)}{A(z)\overline{A(w)}}=
\frac{\phi'(z)K_s(\phi(z),\phi(w))\overline{\phi'(w)}}
{\phi'(z)\cdot(A_s\circ\phi)(z)\ \overline{\phi'(w)\cdot(A_s\circ\phi)(w)}} \\
= \frac{K_s(\phi(z),\phi(w))}
{A_s(\phi(z))\,\overline{A_s(\phi(w))}}.
\endgather$$
If $G_1$ and $G_2$ form a primitive pair for the double of $\O_s$,
then Theorem~1.1 yields that the last function in the string of
equations is a rational combination of $G_1\circ\phi$ and
$G_2\circ\phi$, and this shows that the Bergman kernel $K(z,w)$ is
generated by a primitive pair in the generalized sense that we
defined in \S1.  This completes the proof of Theorem~1.2.

We now turn to the proof of Theorem~1.3.  We have shown that
$\Phi'$ is a linear combination of functions of the
form $K_s^{(m)}(z,W_j)$ where $m\ge0$ and $W_j=\phi(w_j)$.  All
functions of this form extend holomorphically past the boundary of
$\O_s$ (see \cite{4, p.~41,133}).  Hence $\Phi'$ extends
holomorphically past the boundary of $\O_s$ and, consequently, so does
$\Phi$.  It follows that $\O$ is a bounded domain and that the
boundary of $\O$ is piecewise real analytic.
The singular points, if any, in the boundary of $\O$ are given as the
images of boundary points under $\Phi$ where $\Phi'$ has a simple zero.
(Notice that for $\Phi$ to be one-to-one on $\O_s$, $\Phi'$ cannot
have any zeroes of multiplicity greater than one on the boundary.)
It is clear that $\phi$ extends continuously to the boundary of $\O$.
We will complete the proof of Theorem~1.3 during the course of the
next proof when we show that $\Phi$ extends meromorphically to the
double of $\O_s$.
\enddemo

\demo{Proof of Theorem~1.4}
As in the preceding proof, we suppose that $\O$ is an $n$-connected
quadrature domain of finite area in the plane such that no boundary
component is a point and we suppose that $\phi$ is a holomorphic mapping
which maps $\O$ one-to-one onto a bounded domain $\O_s$ bounded by
$n$ smooth real analytic curves.  Again, let $\Phi$ denote $\phi^{-1}$.
We now claim that $\Phi$ extends to the double of $\O_s$ as a
meromorphic function.  To see this, we use the inhomogeneous Cauchy
integral formula,
$$u(z)=\frac{1}{2\pi i}\int_{w\in b\O_s}\frac{u(w)}{w-z}\ dw +
\frac{1}{2\pi i}
\iint_{w\in \O_s}\frac{\dee u/\dee\bar w}{w-z}\ dw\wedge d\bar w$$
with $\overline{\Phi}$ in place of $u$ to obtain
$$
\overline{\Phi(z)}=
\frac{1}{2\pi i}\int_{w\in b\O_s}\frac{\overline{\Phi(w)}}{w-z}\ dw +
\frac{1}{2\pi i}
\iint_{w\in \O_s}\frac{\overline{\Phi'(w)}}{w-z}\ dw\wedge d\bar w.$$
We now let $z$ approach the boundary of $\O_s$ from the inside and
use a method developed in \cite{4} to determine the boundary
values of the two integrals.  Since the boundary of $\O_s$ is real
analytic and since $\Phi$ extends holomorphically past the boundary
of $\O_s$, the Cauchy-Kovalevskaya theorem (see \cite{4, p.~39})
yields that there is a function $v$ which is real analytic in a
neighborhood of the boundary of $\O_s$, which vanishes on the
boundary, and such that $\frac{\dee v}{\dee\bar w}\equiv\overline{\Phi'(w)}$
on a neighborhood of the boundary.  We may extend $v$ to be $C^\infty$
smooth in $\O_s$.  Now apply the inhomogeneous Cauchy Integral Formula
to $\overline{\Phi}-v$ to obtain
$$
\overline{\Phi(z)}-v(z)=
\frac{1}{2\pi i}\int_{w\in b\O_s}\frac{\overline{\Phi(w)}}{w-z}\ dw +
\frac{1}{2\pi i}
\iint_{w\in \O_s}\frac{\overline{\Phi'(w)}-\frac{\dee v}{\dee\bar w}}{w-z}
\ dw\wedge d\bar w.$$
The function given as the integral over the boundary of $\O_s$
is a holomorphic function $H(z)$ on $\O_s$ which extends $C^\infty$
smoothly to the boundary (see \cite{4, p.~7}).
The function given as the integral over $\O_s$ is also $C^\infty$
up to the boundary of $\O_s$ because the integrand has compact
support.  Furthermore, the boundary values of this function agree with the
boundary values of
$$\iint_{w\in \O_s}\frac{\overline{\Phi'(w)}-\frac{\dee v}{\dee\bar w}}{w-z}
\ dw\wedge d\bar w.$$
as $z$ approaches the boundary of $\O_s$ from the {\it outside\/}
of $\O_s$.  It was shown in the course of
proving Theorem~1.3 that $\Phi'$ is equal to a linear combination
of functions of $z$ of the form $K_s^{(m)}(z,w_j)$.  Note that this
integral is a constant times the $L^2$ inner product of
$\Phi'-\frac{\dee\bar v}{\dee w}$, and the holomorphic function
$1/(w-z)$ of $w$ when $z$ is outside of $\O_s$.  Functions of the
form $\dee\bar v/\dee w$ are orthogonal to smooth holomorphic functions
when $v$ vanishes on the boundary via integration by parts.
Hence, the boundary values of the integral over $\O_s$ are given as
a linear combination of terms of the form $1/(w_j-z)^m$,
i.e., a rational function of $z$.  Hence, for $z$ in the boundary
of $\O_s$, we have shown that
$$\overline{\Phi(z)}= H(z)+\Cal R(z)$$
where $H$ is a holomorphic function on $\O_s$ which is $C^\infty$
smooth up to the boundary and $\Cal R(z)$ is a rational function of
$z$.  (In fact, the poles of $R(z)$ only fall at the points $\phi(w_j)$
where $w_j$ are the
points appearing in the quadrature identity for $\O$.)  This last
identity reveals that $\Phi$ extends to the double of $\O_s$ as
a meromorphic function, and this completes the proof of Theorem~1.3,
as promised.

Since $z=\Phi(\phi(z))$ and $\Phi$ extends meromorphically to the
double of $\O_s$, it follows that $z$ extends meromorphically to
the double of $\O$ in the generalized sense that we use in this
paper.

Let the symbol $\Phi$ also denote the function defined on the double
of $\O_s$ which is the meromorphic extension of $\Phi$.  Let $R$ denote
the antiholomorphic reflection function on the double of $\O_s$ which maps
$\O_s$ to reflected copy in the double of $\O_s$ and let
$G=\overline{\Phi\circ R}$.  It is now an easy matter to see that $\Phi$
and $G$ form a primitive pair for the double of $\O_s$.  Indeed, $G$
has poles only at finitely many points in the double which fall in
the $\O_s$ side.  Choose a complex number $w_0$ sufficiently close to
the point at infinity so that the set $\Cal S=G^{-1}(w_0)$ consists
of finitely many points in $\O_s$ such that $G$ takes the value $w_0$
with multiplicity one at each point in $\Cal S$.  Since $\Phi$ is
one-to-one on $\O_s$, it follows that $\Phi$ separates the points of
$G^{-1}(w_0)$, and this implies that $\Phi$ and $G$ form a primitive
pair (see Farkas and Kra \cite{13}).  Consequently, there is an
irreducible polynomial $P(z,w)$ such that $P(\Phi(z),G(z))\equiv0$
on the double of $\O_s$.

The Schwarz function for $\O$ is now given as $S(z)=G(\phi(z))$ and
$z$ and $S(z)$ form a primitive pair for $\O$.  It follows that $S(z)$
is meromorphic on $\O$ and continuous up to the boundary with boundary
values equal to $\bar z$.  By composing the polynomial identity
$P(\Phi(z),G(z))\equiv0$ with $\phi$, we see that $P(z,S(z))\equiv0$
on $\O$, and this shows that $S(z)$ is an algebraic function.
When this identity is restricted to the boundary, we see that
$P(z,\bar z)=0$ when $z$ is in the boundary, i.e., that the boundary
of $\O$ is contained in an algebraic curve.  (Gustafsson refined this
argument to show that the boundary is in fact equal to the algebraic
curve minus perhaps finitely many points.)

We now turn to examine the function $Q(z)$.  Let $p(z)$ denote the
sum of the principal parts of $S(z)$ in $\O$.  Apply the Cauchy
integral formula to $S(z)-p(z)$ and use the fact that $S(z)=\bar z$ on
the boundary to see that $S(z)-p(z)=Q(z)$ plus a linear combination
of integrals of the form $\int_{b\O}\frac{1}{(w-w_j)^k(w-z)}\ dw$
where the $w_j$ are points in $\O$ where $S$ has poles.  But all
such integrals are zero.  Hence $S(z)=p(z)+Q(z)$.  Since $S$ is
algebraic, so is $Q$.  Since $z$ and $S$ generate the meromorphic
functions on the double of $\O$, and since $p(z)$ is rational,
it follows that $z$ and $Q(z)$ also generate the meromorphic
functions on the double of $\O$.  This completes the proof of
Theorem~1.4

\enddemo

\demo{Proof of Theorem~1.5}
Theorem~1.4 together with Theorem~1.2 yield that the Bergman kernel
$K(z,w)$ for $\O$ is a rational combination of $z$, $S(z)$, $\bar w$,
and $\overline{S(w)}$.  Thus, the Bergman kernel is algebraic.
Similarly $K(z,w)$ is a rational combination of $z$ and $Q(z)$.
Since any proper holomorphic mapping of $\O$ onto the unit disc
extends meromorphically to the double, it follows that all such maps
are rational combinations of $z$ and $S(z)$, and consequently, they
are algebraic.

It is proved in \cite{5} that if the Bergman kernel is algebraic,
then so is the Szeg\H o kernel, and so are the classical functions
$F_j'$.
\enddemo

\demo{Proof of Theorem~1.6}
We have seen that the kernels $K(z,w)$ and $S(z,w)^2$ and the proper
holomorphic maps to the unit disc and the functions $F_j$
are all generated by $z$ and $S(z)$, and since these functions
are equal to $z$ and $\bar z$, respectively on the boundary, we may
deduce most of the rest of the claims made in Theorem~1.6.  To
finish the proof, note that identity (2.2) yields that
$$T(z)^2=-\frac{S(a,z)^2}{L(z,a)^2}$$
where $a$ is an arbitrary point chosen and fixed in $\O$.
The function $S(z,a)^2$ is a rational function of $z$ and
$\bar z$ on the boundary.  Identity (2.3) yields that
$L(z,a)^2=S(z,a)^2/f_a(z)^2$, and so $L(z,a)^2$ is also a rational
function of $z$ and $\bar z$ on the boundary.  Finally, it follows
that $T(z)^2$ is a rational function of $z$ and $\bar z$ on the
boundary.
\enddemo

\demo{Proof of Theorem~1.9}
Suppose $\O$ is a finitely connected quadrature domain of finite area.
We know that $\bar z$ is equal to the Schwarz function $S(z)$ on the
boundary and that $S(z)$ is meromorphic on $\O$ with finitely many
poles.  Notice that
$$\gather
\iint_\O z\ \overline{h(z)}\ dA =
\frac{i}{2}
\iint_\O z\ \overline{h(z)}\ dz\wedge d\bar z \\
=
\frac{i}{4}
\iint_\O \frac{\dee}{\dee z}\left(z^2\ \overline{h(z)}\right)
\ dz\wedge d\bar z \\
=
\frac{i}{4}
\int_{b\O} z^2\ \overline{h(z)}\ d\bar z = 
\frac{i}{4}
\int_{b\O} \overline{S(z)^2}\ \overline{h(z)}\ d\bar z,
\endgather$$
and the Residue Theorem yields that this last integral is equal to
a fixed linear combination of values of $h$ and finitely many of its
derivatives at the points in $\O$ where $S(z)$ has poles.  This shows
that the function $z$ is a linear combination of functions of the
form $K^{(m)}(z,w_k)$ where $w_k$ are points where $S(z)$ has poles.
Since the constant function $1$ is also given by a linear combination
as in formula (1.2), we may form a quotient to get a Bergman
representative mapping which is equal to the function $z$.
\enddemo

\subhead 4. Gustafsson's Theorem\endsubhead
Bj\"orn Gustafsson has granted me permission to include here
his proof of his discovery that every Gustafsson map on a smoothly
bounded domain is a Bergman representative map.  After I give
Gustafsson's elegant proof, I will offer an alternative, more
pedestrian proof that sheds additional light on this phenomenon.

The main tool in Gustafsson's argument is the following lemma.

\proclaim{Lemma 4.1}
Suppose that $\O$ is a bounded finitely connected
domain bounded by simple closed $C^\infty$ smooth curves.
If $G$ is a holomorphic function on $\O$ which extends meromorphically
to the double of $\O$ and has no poles on $\Obar$, then $G'$ must be equal
to a complex linear combination of functions of $z$ of the form
$K^{(m)}(z,w_k)$ where $w_k$ are points in $\O$.
\endproclaim

\demo{Proof of Lemma~4.1}
Let $G$ denote both the function on $\O$ and its extension to the double
$\Oh$ of $\O$.  Note that $dG$ is a meromorphic differential on $\Oh$, and
that since it is exact, it is free of residues.  Identity  (3.1) can be
used in a standard way to establish the well know fact that differentials
of the form $K^{(m)}(z,w_k)\,dz$ extend meromorphically to $\Oh$.  Since
the function $\Lambda(z,w)$ has a double pole at $z=w$ with no residue
term, it is possible to choose a linear combination of the extensions of
$K^{(m)}(z,w_k)\,dz$ so that the poles of $dG$ on the back side of $\Oh$
are exactly cancelled by the sum.  (We shall say more about this argument
below.)  Let $\Cal K(z)\,dz$ denote such a linear
combination.  Now $dG - \Cal K\,dz$ is a holomorphic differential on $\Oh$.
All the periods of $dG$ are zero, and it is well known that the
$\beta$-periods of differentials of the form $\Cal K\,dz$ vanish, i.e.,
the periods that go across $\O$ from one boundary curve to another and
then return along the backside along the reflected curve.  (This fact follows
from the relationships between the Bergman kernel and the $\Lambda$ kernel
and the Green's function.  It is explained in Schiffer and Spencer
\cite{16, p.~101-105}.  Also, the arguments can be found in Gustafsson's
paper \cite{14} in the proofs of Theorems~1 and~2.)

Now a holomorphic differential vanishes if all its $\beta$-periods, or
if all its $\alpha$-periods, vanish.  Hence $dG-\Cal K\,dz$ is zero and
it follows that $G'=\Cal K$ on $\O$.  This completes the proof.
\enddemo

I call a function on a bounded finitely connected domain bounded by
smooth curves a Gustafsson mapping if it is holomorphic
and one-to-one on the domain and extends meromorphically to the
double and has no poles on the boundary of the domain.  Gustafsson
mappings effect a conformal change of variables from the given domain
to a quadrature domain.

\proclaim{Theorem 4.2}
Suppose that $\O$ is a bounded finitely connected
domain bounded by simple closed $C^\infty$ smooth curves.
If $g$ is a Gustafsson function on $\O$, then $g$ is equal
to a Bergman representative mapping.
\endproclaim

\demo{Proof of Theorem 4.2}
If $g$ is a Gustafsson map, apply Lemma~4.1 to $g$ and to $\frac12 g^2$
to get $g'=\Cal K_2$ and $g g'=\Cal K_1$.  Since $g$ is one-to-one, the
equality $g'=\Cal K_2$ shows that the linear combination $\Cal K_2$ is
non-vanishing on $\Obar$.  Now $g$ is given by the quotient
$\Cal K_1/\Cal K_2$, which is a Bergman representative mapping.
\enddemo

I shall now give an alternative proof of Lemma~4.1 that sheds further
light on the property proved in \cite{9} that, on a smooth quadrature domain, 
if $G$ is a function on the domain that extends meromorphically to the
double, then $G'$ is also a function on the domain that extends
meromorphically to the double.  Suppose that $\O$ is a bounded finitely
connected domain bounded by simple closed $C^\infty$ smooth curves,
and suppose that $G$ is a holomorphic function on $\O$ that extends
meromorphically to the double of $\O$ and has no poles on $\Obar$.
In this case, there is a meromorphic function $H$ on $\O$ which extends
smoothly up to the boundary such that $G=\overline{H}$ on $b\O$.
Let $z(t)$ denote a parameterization of a boundary curve of $\O$.
Since $G(z(t))=\overline{H(z(t))}$, we may differentiate with respect
to $t$ and divide the result by $|z'(t)|$ to obtain that
$G'(z)T(z)=\overline{H'(z)T(z)}$ for $z\in b\O$.
This last identity is very similar to identity (3.1).  Indeed,
we may differentiate (3.1) with respect to $w$ and rewrite it
to obtain
$$K^{(m)}(z,w)T(z)=
-\overline{\Lambda^{(m)}(z,w)T(z)}$$
for $z\in b\O$, where
$\Lambda^{(m)}(z,w)$ denotes the $m$-th derivative of
$\Lambda(z,w)$ with respect to $w$ (and of course
$\Lambda^{(0)}(z,w)=\Lambda(z,w)$).  The singular part of
$\Lambda(z,w)$ is a constant times $(z-w)^{-2}$.  Since
$H'$ is the derivative of a meromorphic function, the
poles of $H'$ are double or more.  Hence,
there is a unique linear combination $\Cal L$ of the functions
$\Lambda^{(m)}(z,w)$ so that the principal parts at the poles
of $\Cal L$  agree with the principal parts of $H'$ at each
pole in $\O$.  If
$\Cal L(z)=\sum_{j=1}^N\sum_{m=1}^{n_j}c_{jm}\Lambda^{(m)}(z,w_j)$,
let
$\Cal K(z)=-\sum_{j=1}^N\sum_{m=1}^{n_j}\bar c_{jm}K^{(m)}(z,w_j)$.
Notice that $\Cal K(z)T(z)= \overline{\Cal L T(z)}$ on $b\O$.
Now $(G'-\Cal K)T=\overline{(H'-\Cal L)T}$ on $b\O$ where both
$G'-\Cal K$ and $H'-\Cal L$ are holomorphic on $\O$ and extend
smoothly to the boundary.  This implies
$(G'-\Cal K)T$ is both orthogonal to the Hardy space and conjugates
of functions in the Hardy space.  Consequently, a theorem of Schiffer
yields that
$G'-\Cal K=\sum_{j=1}^{n-1}c_jF_j'$ for some constants $c_j$.   (See
\cite{4, p.~80} for a proof of this result that proves, rather than
assumes, that the zeroes of the Szeg\"o kernel are simple zeroes.)
Now $\Cal K(z)$ is $\dee/\dee z$ of a linear combination $\Cal G$ of
functions of the form $(\dee^m/\dee\bar w^m)G(z,w)$ where $G(z,w)$ is the
classical Green's function.  All linear combinations of this form
vanish on the boundary in the $z$ variable.  Also,
$F_j'=2(\dee/\dee z)\omega_j$.  Hence
$G'-\Cal K=\sum_{j=1}^{n-1}c_jF_j'$ yields that
$G-\Cal G-2\sum_{j=1}^{n-1}c_j\omega_j$ is antimeromorphic on $\O$.  But
$G=\overline{H}$ on the boundary.  Hence, the boundary values of
$\Cal G+2\sum_{j=1}^{n-1}c_j\omega_j$ agree with the boundary
values of a function which is antimeromorphic on $\O$ and which extends
smoothly to the boundary.  Since this meromorphic function vanishes on
one of the boundary curves, it must be identically zero.  This forces us
to conclude that all the $c_j$'s are zero, and we have proved
Gustafsson's theorem that $G'=\Cal K$.

It should be remarked that the argument just given can be run in reverse
to yield a converse to Gustafsson's lemma.  Indeed, if $\Cal K$ is a
linear combination of the form used above that has vanishing periods along
all the boundary curves of $\O$, then an analytic antiderivative of
$\Cal K$ must extend meromorphically to the double.

Gustafsson's lemma can also be routinely generalized to yield the
following result.  Suppose that $\O$ is a bounded finitely connected
domain bounded by simple closed $C^\infty$ smooth curves.  If $G$ is
a {\it meromorphic\/} function on $\O$ which extends meromorphically
to the double of $\O$, then $G'(z)$ must be equal to a function of
the form $$\sum_{k=1}^N\sum_{m=1}^{p_k} a_{km} K^{(m)}(z,x_k) +
\sum_{k=1}^M\sum_{m=1}^{q_k} b_{km}\Lambda^{(m)}(z,y_k),$$
where the $x_k$ are points in $\Obar$ and the $y_k$ are points in $\O$.
The points $x_k$ are used to cancel the poles of the extension of $G'$
on the back side of the double and the $y_k$ are used to cancel the poles of
$G'$ on $\O$.  These results allow the field of meromorphic functions
on the double to be handled like a linear space in many instances.

In case the domain under study is a quadrature domain of finite area
with smooth boundary, then we know that the Schwarz function $S(z)$
agrees with $\bar z$ on the boundary.  We can differentiate the identity
$z=\overline{S(z)}$ along the boundary as we did above to see that
$T(z)=\overline{S'(z)T(z)}$ for $z$ in $b\O$.  Notice that this identity
reveals that $|S'(z)|=1$ on the boundary.  Now, if $G$ is a meromorphic
function on $\O$ which extends meromorphically to the double, we may
write $G'T=\overline{H'T}$ on the boundary as we did above and divide
by the identity for the Schwarz function to obtain that $G'(z)$ is
equal to the conjugate of $H'(z)/S'(z)$ for $z$ in the boundary.  This
yields another way to see that, on a smooth quadrature domain, any
meromorphic function $G$ on the domain that extends meromorphically
to the double has the property that its derivative $G'$ also extends
meromorphically to the double.  This proof also has the virtue that
it gives an explicit formula for the extension of $G'$ to the double.

\enddocument